\documentclass[a4paper,12pt]{amsart}
\usepackage{amscd,amssymb}
\author{D.~A.~Timashev}
\address{
Department of Higher Algebra\\
Faculty of Mechanics and Mathematics\\
Moscow State University\\
119992 Moscow, Russia}
\email{timashev@mech.math.msu.su}
\urladdr{http://mech.math.msu.su/department/algebra/staff/timashev}
\thanks{Supported by CRDF grant RM1--2543--MO--03 and RFBR grant 05--01--00988.}
\title[Equivariant symplectic geometry of cotangent bundles]
{Equivariant symplectic geometry\\ of cotangent bundles, II}
\date{February 14, 2005}

\newcommand{\ab}{\mathfrak a}
\newcommand{\g}{\mathfrak g}
\newcommand{\h}{\mathfrak h}
\newcommand{\lv}{\mathfrak l}
\newcommand{\m}{\mathfrak m}
\newcommand{\M}{\mathfrak M}
\newcommand{\p}{\mathfrak p}
\newcommand{\q}{\mathfrak q}
\newcommand{\CC}{\mathbb C}
\newcommand{\PP}{\mathbb P}
\newcommand{\N}{\mathcal N}
\newcommand{\U}{\mathcal U}
\newcommand{\X}{X_0}
\newcommand{\Xx}{X_1}
\newcommand{\Z}{Z_0}
\newcommand{\Zz}{Z_1}
\let\phi\varphi

\newcommand{\ad}{\mathop{\mathrm{ad}}}
\newcommand{\codim}{\mathop{\mathrm{codim}}}
\newcommand{\conn}[1]{#1^{\circ}}
\newcommand{\diag}{\mathop{\mathrm{diag}}}
\renewcommand{\Im}{\mathop{\mathrm{Im}}}
\newcommand{\Ker}{\mathop{\mathrm{Ker}}}
\newcommand{\nil}[1]{#1_{\mathrm{n}}}
\newcommand{\sms}[1]{#1_{\mathrm{s}}}
\newcommand{\pr}{\mathrm{pr}}
\newcommand{\Ru}[1]{#1_{\mathrm{u}}}
\newcommand{\itimes}[1]{*_{#1}}
\newcommand{\ftimes}[1]{\mathop{\underset{#1}{\times}}}

\newtheorem{lemma}{Lemma}
\newtheorem{proposition}{Proposition}
\newtheorem{theorem}{Theorem}
\newtheorem{corollary}{Corollary}
\newtheorem*{Corollary}{Corollary}
\newtheorem{claim}{Claim}
\newtheorem*{Claim}{Claim}
\theoremstyle{definition}
\newtheorem{definition}{Definition}
\newtheorem{example}{Example}
\theoremstyle{remark}
\newtheorem{remark}{Remark}

\newcommand{\ND}{\textup{(ND)}}
\newcommand{\DX}{\textup{(D${}_X$)}}
\newcommand{\DY}{\textup{(D${}_Y$)}}

\begin{document}

\begin{abstract}
We examine the structure of the cotangent bundle $T^{*}X$ of an
algebraic variety $X$ acted on by a reductive group $G$ from the
viewpoint of equivariant symplectic geometry. In particular, we
construct an equivariant symplectic covering of $T^{*}X$ by the
cotangent bundle of a certain variety of horospheres in~$X$, and
integrate the invariant collective motion on~$T^{*}X$. These
results are based on a ``local structure theorem'' describing
the action of a certain parabolic in $G$ on an open subset
of~$X$, which is interesting by itself.
\end{abstract}

\maketitle

\section*{Introduction}

An important class of symplectic manifolds with a Hamiltonian
group action is formed by cotangent bundles of manifolds acted
on by Lie groups. Cotangent bundles arise as phase spaces for
many important Hamiltonian dynamical systems with symmetries.
Therefore it is an important problem to study the equivariant
symplectic geometry of cotangent bundles. In particular, one may
address the problem of constructing a symplectic manifold which
is locally equivariantly isomorphic to the cotangent bundle but
has a simpler structure than the latter one. Morally, this should
help in integrating Hamiltonian dynamical systems on cotangent
bundles.

In this paper, we study the equivariant geometry of cotangent
bundles in the framework of algebraic geometry. More precisely,
let $X$ be a smooth algebraic variety over complex numbers (or,
more generally, over any algebraically closed field of
characteristic zero) acted on by a connected reductive
group~$G$. We examine the natural $G$-action on the cotangent
bundle $T^{*}X$ from the viewpoint of transformation groups and
symplectic geometry.

This problem has attracted the attention of several researchers.
In particular, F.~Knop studied the moment map
$\Phi:T^{*}X\to\g^{*}$ in~\cite{W_X}. Using the moment map, he
obtained deep results on the geometry of the action $G:T^{*}X$
including: the existence and the description of the stabilizer
in general position; the description of the collective functions
(which are the integrals for any $G$-invariant Hamiltonian system
on $T^{*}X$); the relation between symplectic invariants of
$T^{*}X$ (corank, defect) and important invariants of $X$
(complexity, rank), which play a significant role, e.g., in
studying equivariant embeddings of $X$; etc. In~\cite{inv.mot}
Knop studied the invariant collective motion on $T^{*}X$, i.e.,
the flow generated by the skew-gradients of invariant collective
functions, under restriction that $X$ is ``non-degenerate''.
(This class includes, e.g., all quasiaffine varieties.)

On the other hand, E.~B.~Vinberg \cite{T*X} constructed an
equivariant symplectic rational Galois covering of $T^{*}X$ by
the cotangent bundle of the variety of generic horospheres
provided that $X$ is quasiaffine. (A horosphere in $X$ is an
orbit of a maximal unipotent subgroup of~$G$.) Since the variety
of horospheres and its cotangent bundle have a relatively simple
structure from the point of view of $G$-action, this result
solves the problem posed above.

Our main objective here is to generalize Vinberg's construction to
arbitrary~$X$. The general scheme of reasoning is the same as in
\cite{T*X}, so that our paper may be regarded as a direct
continuation of \cite{T*X}.%
\footnote{The author is grateful to E.~B.~Vinberg for his kind
permission to entitle this paper as Part~II of~\cite{T*X}.}
On the other hand, the main technical tool in the paper is a
refined version of the so-called ``local structure theorem'' (see
Section~\ref{loc.struct}), which in its turn yields simple proofs
and generalizations of some results from \cite{W_X},
\cite{inv.mot}. Roughly speaking, the main results of \cite{W_X}
stem from a ``quantization'' argument, i.e., by passing from
$T^{*}X$ to differential operators on~$X$. We reprove them by pure
``classical'' arguments in the spirit of \cite{inv.mot} but
dropping the ``non-degeneracy'' assumption. We are also able to
deduce some results of Knop on invariant collective motion in full
generality.

The structure of the paper is as follows. In
Section~\ref{polar} we describe a general construction from
symplectic geometry used by Vinberg to construct his rational
covering. Then we recall Vinberg's result and explain why it
does not generalize directly. In Section~\ref{loc.struct} we
prove the refined local structure theorem. Using this theorem,
we examine the geometry of $T^{*}X$ in Section~\ref{cotan}. In
particular, we describe the image of the moment map and the
stabilizer in general position. An equivariant symplectic
rational Galois covering of $T^{*}X$ by the cotangent bundle of
a certain variety of ``degenerate'' horospheres is constructed
in Section~\ref{horosph}. In Section~\ref{inv.mot} we study the
invariant collective motion on $T^{*}X$ and show how it can be
integrated after lifting via the above covering.

\section*{Notation and conventions}

All algebraic varieties and groups are defined over one and the
same algebraically closed base field of characteristic zero.

If $X$ is a smooth variety, then $TX$ (resp.~$T^{*}X$) is its
(co)tangent bundle and $NY$ (resp.~$N^{*}Y$) is the (co)normal
bundle of a smooth subvariety $Y\subseteq X$. The fibres at a
point $x$ are denoted by the respective subscript: $T_xX$,
$N^{*}_xY$, etc.

Given a vector space~$V$, the dual space is denoted by~$V^{*}$,
and $\langle\omega,v\rangle$ is the pairing of $v\in V$ and
$\omega\in V^{*}$.

If $\phi:X\to Y$ is a map and $f$ is a function on~$Y$, then
$\phi^{*}f$ denotes its pullback to~$X$:
$\phi^{*}f(x)=f(\phi(x))$.

Algebraic groups are denoted by capital Latin letters and their
Lie algebras are denoted by the respective small Gothic letters.
The identity component of an algebraic group $G$ is denoted
by~$\conn{G}$. For any subset $S$ in $G$ or in $\g$ let $Z(S)$
be the centralizer and $N(S)$ the normalizer of $S$ in~$G$. The
additive Jordan decomposition of $\xi\in\g$ is written as
$\xi=\sms\xi+\nil\xi$ with $\sms\xi$ semisimple and $\nil\xi$
nilpotent.

An action of a group $G$ on a set $X$ is denoted by ${G:X}$. Let
$X^G$ denote the set of fixed points, $Gx$~the orbit of $x\in
X$, and $G_x$ the stabilizer of $x$ in~$G$. For an algebraic
group action, $\xi{x}$~is the velocity vector of $\xi\in\g$ at
$x\in X$, and $\g{x}=T_xGx$ is the tangent space to the orbit.

When we speak of an action $G:\g$ (or $G:\g^{*}$), it is always
assumed that the action is (co)adjoint.

If $H\subseteq G$ is a subgroup and $Y$ is an $H$-set, then
$G\itimes{H}Y$ denotes the homogeneous fibre bundle over $G/H$
with the fibre $Y$ over the base point. It can be defined as
$(G\times Y)/H$, where $H:G\times Y$ is defined by
$h(g,y)=(gh^{-1},hy)$. Let $g*y$ denote the point in
$G\itimes{H}Y$ corresponding to~$(g,y)$. The construction of a
homogeneous fibre bundle is well defined in the category of
algebraic varieties.

\section{Polarization of cotangent bundles}\label{polar}

\subsection{}\label{skew}

In \cite[\S4]{T*X} Vinberg introduced a general construction
from symplectic geometry which relates the cotangent bundles of
two different varieties. The desired symplectic covering of a
cotangent bundle is its particular application. So for
convenience of the reader we recall this construction here.

Let $X$, $Y$, and $Z\subseteq X\times Y$ be smooth irreducible
algebraic varieties. Assume that the projections $p:Z\to X$,
$q:Z\to Y$ are smooth surjective maps. Then $Z$ can be regarded
as a family of smooth subvarieties $Z_y=\{x\mid(x,y)\in Z\}
\subseteq X$ of equal dimension parametrized by points $y\in Y$,
or similarly, as a family of subvarieties $Z_x\subseteq Y$.

For any $z=(x,y)\in X\times Y$ there are canonical isomorphisms
\begin{equation*}
T_z(X\times Y)\simeq T_xX\oplus T_yY,\qquad
T^{*}_z(X\times Y)\simeq T^{*}_xX\oplus T^{*}_yY.
\end{equation*}
For any $\alpha\in T^{*}_z(X\times Y)$ we denote by
$\alpha',\alpha''$ its projections to $T^{*}_xX$ and $T^{*}_yY$,
respectively.

\begin{definition}
The \emph{skew conormal bundle} of $Z$ is
\begin{equation*}
SN^{*}Z=\{\alpha\in T^{*}_z(X\times Y)\mid z\in Z,\
\alpha'-\alpha''=0\text{ on }T_zZ\}.
\end{equation*}
\end{definition}

\begin{remark}
The skew conormal bundle is obtained from the usual conormal
bundle $N^{*}Z$ by an automorphism of $T^{*}(X\times Y)\simeq
T^{*}X\times T^{*}Y$, namely, by multiplying the covectors over
$Y$ by~$-1$. However it is more convenient to consider the skew
conormal bundle instead of the conormal bundle in order to avoid
superfluous signs in some formul{\ae}.
\end{remark}

There is a commutative diagram
\begin{equation*}
\begin{CD}
T^{*}X @<{\hat{p}}<< SN^{*}Z @>{\hat{q}}>> T^{*}Y \\
 @VVV                 @VVV                  @VVV  \\
   X      @<{p}<<       Z       @>{q}>>       Y,
\end{CD}
\end{equation*}
where $\hat{p},\hat{q}$ take $\alpha$ to $\alpha',\alpha''$,
respectively.

\begin{remark}
Consider the disjoint union $\bigsqcup_{y\in Y}N^{*}Z_y$ of the
conormal bundles of the subvarieties $Z_y$ in~$X$. It has a
natural structure of a subbundle in the pullback of $T^{*}X$
to~$Z$. By definition, this bundle consists of pairs
$(\alpha',y)$ such that $\alpha'\in T^{*}_xX$, $Z_y\ni x$, and
$Z_y$ is tangent to $\Ker\alpha'$ at~$x$. We call such pairs the
\emph{polarized covectors} over~$X$, and the whole bundle is
called the \emph{polarized cotangent bundle} of $X$ with respect
to the family of subvarieties~$Z_y$.

It is easy to see that $\hat{p}$ induces an isomorphism
\begin{equation}\label{skew=pol}
SN^{*}Z\simeq\bigsqcup_{y\in Y}N^{*}Z_y.
\end{equation}
Indeed, for any $z=(x,y)\in Z$ the linear map
$\hat{p}:SN^{*}_zZ\to N^{*}_xZ_y$ is isomorphic, because
$T_zZ\cap T_z({X\times\{y\}})=T_z(Z_y\times\{y\})$ and
$T_zZ+T_z(X\times\{y\})=T_z(X\times Y)$ by smoothness and
surjectivity of~$q$. Similarly,
\begin{equation*}
SN^{*}Z\simeq\bigsqcup_{x\in X}N^{*}Z_x.
\end{equation*}
\end{remark}

Recall that the canonical symplectic structure on $T^{*}X$
arises from a certain 1-form $\ell'$, called the \emph{action
form}. Given $\alpha'\in T^{*}_xX$ and $\nu'\in
T_{\alpha'}T^{*}X$, we put
$\ell'(\nu')=\langle\alpha',\xi'\rangle$, where $\xi'\in T_xX$
is the projection of~$\nu'$. Then the symplectic form on
$T^{*}X$ is defined as $\omega'=d\ell'$. Similarly, one defines
the action form $\ell''$ and the symplectic form $\omega''$
on~$T^{*}Y$.

The action forms $\ell',\ell''$ lift to one and the same 1-form
$\ell$ on~$SN^{*}Z$. Indeed, given any $\alpha\in SN^{*}Z$,
$z=(x,y)$, and any $\nu\in T_{\alpha}SN^{*}Z$, denote
$\nu'=d\hat{p}(\nu)\in T_{\alpha'}T^{*}X$,
$\nu''=d\hat{q}(\nu)\in T_{\alpha''}T^{*}Y$, and by
$\xi,\xi',\xi''$ the projections of $\nu,\nu',\nu''$ to
$T_zZ,T_xX,T_yY$, respectively. Then
$\langle{\alpha'-\alpha''},\xi\rangle=
\langle\alpha',\xi'\rangle-\langle\alpha'',\xi''\rangle=0$
implies $\ell'(\nu')=\ell''(\nu'')$.

It follows that $\omega',\omega''$ lift to one and the same
closed (but maybe degenerate) 2-form $\omega=d\ell$
on~$SN^{*}Z$.

\begin{lemma}\label{symp}
If $\dim X=\dim Y$, then the following conditions are equivalent:
\begin{enumerate}
\item[\ND] $\omega$~is non-degenerate at points in general
position;
\item[\DX] $\hat{p}:SN^{*}Z\to T^{*}X$ is dominant;
\item[\DY] $\hat{q}:SN^{*}Z\to T^{*}Y$ is dominant.
\end{enumerate}
\end{lemma}

\begin{proof}
We have $\dim SN^{*}Z=\dim T^{*}X=\dim T^{*}Y=2\dim X$. If
$\omega$ is generically non-degenerate, then $\hat{p}$ has
finite generic fibres, hence it is dominant by dimension count.
Conversely, a dominant map between varieties of equal dimension
is generically {\'e}tale whence \mbox{\ND$\iff$\DX}. Similarly,
\ND$\iff$\DY.
\end{proof}

Under the conditions of the lemma, $\omega$~defines a symplectic
structure on an open subset of~$SN^{*}Z$, so that
$\hat{p},\hat{q}$ are symplectic rational coverings.

\subsection{}\label{hor}

Assume now that $X$ is equipped with an action of a reductive
connected group~$G$. The induced action $G:T^{*}X$ is
\emph{Hamiltonian}, i.e., it preserves the symplectic structure
and there exists a $G$-equivariant \emph{moment map}
$\Phi:T^{*}X\to\g^{*}$ with the following property: for any
$\xi\in\g$, regarded as a linear function on~$\g^{*}$, the skew
gradient of $\Phi^{*}\xi$ equals the velocity field of $\xi$
on~$T^{*}X$. The moment map is defined by the formula
\begin{equation*}
\langle\Phi(\alpha),\xi\rangle=\langle\alpha,\xi{x}\rangle,
\qquad\forall x\in X,\ \alpha\in T^{*}_xX,\ \xi\in\g.
\end{equation*}
For instance, if $X=G/H$, then
\begin{equation}\label{T*(G/H)}
T^{*}X\simeq G\itimes{H}(\g/\h)^{*}\simeq G\itimes{H}\h^{\perp},
\end{equation}
where $\h^{\perp}$ is the annihilator of $\h$ in~$\g^{*}$, and
the moment map amounts to the coadjoint action:
$\Phi(g*\alpha)=g\alpha$.

A \emph{horosphere} in $X$ is an orbit of any maximal unipotent
subgroup of~$G$. This terminology goes back to I.~M.~Gelfand and
M.~I.~Graev \cite{hor}, and it is justified by an observation that for
$X=S^n(\CC)$ (the $n$-dimensional complex sphere),
$G=SO_{n+1}(\CC)$, the (generic) horospheres are nothing else
but the complexifications of usual horospheres in the
Lobachevsky space~$L^n$.

The set of generic horospheres can be equipped with a structure
of an algebraic $G$-variety of the same dimension as~$X$, see
\ref{BLV.loc.str} for details. (Since the notion of a
``generic horosphere'' is not quite well defined, this variety is
determined only up to a birational equivalence, but this
suffices for our purposes.)

In the notation of~\ref{skew}, let $Y$ be the variety of generic
horospheres in~$X$, and $Z=\{(x,H)\mid x\in H\}\subset X\times
Y$ the incidence variety. (Perhaps, one has to shrink $X,Y,Z$ a
little bit in order to fulfil all necessary requirements on
smoothness.) Applying the construction of \ref{skew} to this
setting, Vinberg proved the following

\begin{theorem}[\cite{T*X}]\label{cover(qaff)}
If $X$ is quasiaffine, then there exists a $G$-equi\-vari\-ant
symplectic rational Galois covering $T^{*}Y\dasharrow T^{*}X$.
\end{theorem}

Actually, Knop proved implicitly in \cite{inv.mot} that under
the assumptions of the theorem $\hat{p}:SN^{*}Z\to T^{*}X$ is a
rational Galois covering, and Vinberg showed that
$\hat{q}:SN^{*}Z\to T^{*}Y$ is birational. The desired covering
is then defined as $\hat{p}\hat{q}^{-1}$.

The variety of horospheres $Y$ and its cotangent bundle $T^{*}Y$
have a simple structure as $G$-varieties compared with that of
$X$ and~$T^{*}X$ (see, e.g., \cite{W_X}, \ref{BLV.loc.str}). Thus
Theorem~\ref{cover(qaff)} gives a good approximation to the
structure of $T^{*}X$ as a symplectic $G$-variety.

However Theorem~\ref{cover(qaff)} does not generalize
na{\"\i}vely to arbitrary~$X$. Indeed, in view of {\DX} and
(\ref{skew=pol}) a necessary condition is that $\bigcup_{H\in
Y}N^{*}H$ be dense in~$T^{*}X$, i.e., any covector in general
position must vanish along the tangent space of a suitable
horosphere. But there are simple counterexamples in the
non-quasiaffine case:

\begin{example}\label{flag}
Let $X=G/P$ be a generalized flag variety (i.e., $P$ is a
parabolic subgroup of~$G$). Generic horospheres are just the
open Schubert cells with respect to various choices of a maximal
unipotent subgroup of~$G$. For instance, $X=\PP^n$,
$G=GL_{n+1}$, and generic horospheres are complements to
hyperplanes. Thus $\bigcup_{H\in Y}N^{*}H$ is the zero bundle
over $X$.

However, if we consider the ``most degenerate'' horospheres,
which are just points in~$X$, then everything becomes fine:
conormal bundles are just cotangent spaces at points of~$X$,
$Y=X$, and the covering $T^{*}Y\to T^{*}X$ is the identity map.
\end{example}

This example suggests a remedy in the general case: to take
for $Y$ a certain variety of non-generic horospheres. This idea
is developed in Section~\ref{horosph} and leads to a
generalization of Theorem~\ref{cover(qaff)}.

\section{Local structure theorem}\label{loc.struct}

\subsection{}\label{BLV.loc.str}

Let $G$ be a connected reductive group acting on an irreducible
algebraic variety~$X$. In this section, we describe the action
of a certain parabolic subgroup of $G$ on an open subset of~$X$.
Results of this kind are called ``local structure theorems'' and
are ubiquitous in the study of reductive group actions. They
arise from Brion, Luna, Vust \cite{BLV}, and Grosshans \cite{Gr},
cf.\ \cite{W_X}, \cite{inv.mot}, \cite{T*X}.

Fix a Borel subgroup $B\subseteq G$ with the unipotent part
$U\subset B$. Let $P\supseteq B$ be the largest subgroup of $G$
which stabilizes all $B$-orbits in general position in~$X$. Consider
a Levi decomposition $P=\Ru{P}\leftthreetimes L$, where
$\Ru{P}\subseteq U$ is the unipotent radical of $P$ and $L$ a
Levi subgroup.

\begin{theorem}[{\cite[2.3]{W_X}, \cite[\S2]{inv.mot}}]
\label{gen.loc.str}
There is an open $P$-stable subset $\X\subseteq X$ and a closed
$L$-stable subset $\Z\subseteq\X$ such that the natural $P$-equivariant map
\begin{equation*}
P\itimes{L}\Z\longrightarrow\X,\qquad p*z\mapsto pz,
\end{equation*}
is an isomorphism.
Furthermore, the kernel
$L_0$ of the action $L:Z$ contains $[L,L]$, the torus
$A=L/L_0$ acts on $\Z$ freely, and $\Z\simeq A\times C$ for a
certain closed subvariety $C\subseteq\Z$, so that
\begin{equation*}
\X\simeq\Ru{P}\times A\times C.
\end{equation*}
\end{theorem}

From this theorem, one deduces that $P$ is exactly the stabilizer
in $G$ of a generic $B$-orbit and $P_0=\Ru{P}\leftthreetimes L_0$
is the stabilizer of a generic $U$-orbit (these orbits are
parametrized by points of the cross-sections $C$ and~$\Z$,
respectively) \cite[\S2]{val}, \cite[\S3]{T*X}, cf.~\ref{deg.hor}.
Moreover, an element $g\in G$ translates a generic $B$- or
$U$-orbit to another generic $B$- or $U$-orbit iff $g\in P$.

Since all maximal unipotent subgroups of $G$ are conjugate, each
horosphere is a $G$-translate of a $U$-orbit. By the above, the
set of generic horospheres, i.e., $G$-translates of $U$-orbits
in~$\X$, is isomorphic to $G\itimes{P}\Z\simeq G/P_0\times C$ as
a $G$-set and inherits the structure of an algebraic $G$-variety
from the latter one. Note that its dimension equals $\dim X$.

\subsection{}

In order to formulate our version of the local structure
theorem, we have to introduce more notation.

We fix a $G$-invariant inner product on $\g$ and identify
$\g^{*}$ with $\g$ by means of this product whenever it is
convenient for our purposes.

The torus $A$ is not a subgroup of~$L$, but its Lie algebra
$\ab$ can be embedded in $\lv$ as the orthocomplement to~$\lv_0$.
Then $M=Z(\ab)$ is a Levi subgroup of $G$ containing~$L$, and
$Q=BM$ is a parabolic subgroup containing~$P$, with a Levi
decomposition $Q=\Ru{Q}\leftthreetimes M$. Also put
$M_0=L_0[M,M]$ and $Q_0=\Ru{Q}\leftthreetimes M_0$; then
$A\simeq M/M_0\simeq Q/Q_0$.

By the superscript ``${}^{-}$'' we indicate opposite parabolic
subgroups (i.e., the parabolics intersecting given parabolics
in specified Levi subgroups) and related subgroups (e.g., their
unipotent radicals). The correlation between the Lie algebras of
the groups introduced above is represented at the picture (the
blocks indicate direct summands of~$\g$):

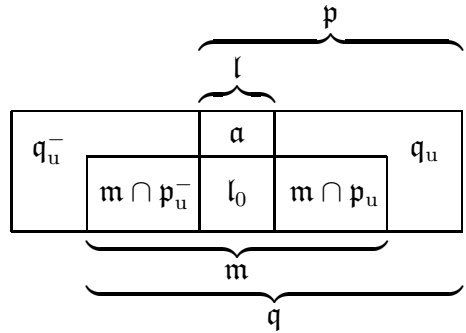
\begin{center}
\unitlength 0.55ex
\linethickness{0.4pt}
\begin{picture}(60.00,43.00)
\put(30.00,18.00){\makebox(0,0)[cc]{$\strut\smash{\lv_0}$}}
\put(30.00,26.00){\makebox(0,0)[cc]{$\ab$}}
\put(43.00,18.00){\makebox(0,0)[cc]{$\strut\smash{\m\cap\Ru\p}$}}
\put(18.00,18.00){\makebox(0,0)[cc]{$\strut\smash{\m\cap\Ru\p^{-}}$}}
\put(55.00,24.00){\makebox(0,0)[cc]{$\strut\smash{\Ru\q}$}}
\put(5.00,24.00){\makebox(0,0)[cc]{$\strut\smash{\Ru\q^{-}}$}}
\put(10.00,13.00){\line(0,1){10.00}}
\put(10.00,23.00){\line(1,0){40.00}}
\put(50.00,23.00){\line(0,-1){10.00}}
\put(60.00,13.00){\line(-1,0){60.00}}
\put(0.00,13.00){\line(0,1){16.00}}
\put(0.00,29.00){\line(1,0){60.00}}
\put(60.00,29.00){\line(0,-1){16.00}}
\put(25.00,13.00){\line(0,1){16.00}}
\put(35.00,13.00){\line(0,1){16.00}}
\put(25.00,30.00){\makebox(0,0)[lb]{$\overbrace{\rule{5.5ex}{0pt}}^{\textstyle\lv}$}}
\put(25.00,37.00){\makebox(0,0)[lb]{$\overbrace{\rule{19.25ex}{0pt}}^{\textstyle\p}$}}
\put(10.00,12.00){\makebox(0,0)[lt]{$\underbrace{\rule{22ex}{0pt}}_{\textstyle\m}$}}
\put(10.00,6.00){\makebox(0,0)[lt]{$\underbrace{\rule{27.5ex}{0pt}}_{\textstyle\q}$}}
\end{picture}
\end{center}

Under the identification $\g\simeq\g^{*}$, the spaces
$\lv,\lv_0,\m,\m_0,\ab$ are self-dual, whereas
$\Ru\p^{-}=\p^{\perp}\simeq(\g/\p)^{*}\simeq\Ru\p^{*}$ and
$\Ru\q^{-}=\q^{\perp}\simeq(\g/\q)^{*}\simeq\Ru\q^{*}$.

\begin{definition}\label{pr}
The \emph{principal stratum} of $\ab$ is
\begin{equation*}
\ab^{\pr}=\{\lambda\in\ab\mid Z(\lambda)=M,\ g\lambda\notin\ab,\
\forall g\in N(M)\setminus N(\ab)\}
\end{equation*}
This is an open subset in $\ab$ complementary to finitely many
linear subspaces (the kernels of the nonzero weights of $\ad\ab$
in $\g$ and the proper subspaces of the form $\ab\cap g\ab$,
$g\in N(M)$).
\end{definition}

Here comes the refined version of the local structure theorem:

\begin{theorem}\label{ref.loc.str}
There is an open $Q$-stable subset $\Xx\subseteq X$ and a closed
$M$-stable subset $\Zz\subseteq\Xx$ such that the natural $Q$-equivariant map
\begin{equation*}
Q\itimes{M}\Zz\longrightarrow\Xx
\end{equation*}
is an isomorphism. Furthermore,
\begin{equation*}
\Zz\simeq M/(M\cap P_0^{-})\times C\simeq M_0/(M_0\cap
P_0^{-})\times A\times C
\end{equation*}
as an $M$-variety and
\begin{equation}\label{X_1}
\Xx\simeq\Ru{Q}\times(M/M\cap P_0^{-})\times C.
\end{equation}
(In the product decompositions, it is always assumed that the
groups act trivially on~$C$.)
\end{theorem}

In comparison with Theorem~\ref{gen.loc.str}, this theorem
displays the local action of a larger parabolic subgroup~$Q$, so
that the locally free action of the complementary part
$M\cap\Ru{P}$ to $\Ru{Q}$ in $\Ru{P}$ is extended to the action of
$M_0$ on a generalized flag variety $M_0/(M_0\cap P^{-})$.

For quasiaffine $X$ one verifies that $M=L$ and $Q=P$
\cite[3.1]{inv.mot}, \cite[\S1]{T*X}, so that
Theorem~\ref{ref.loc.str} specializes to
Theorem~\ref{gen.loc.str}. More generally, this property
characterizes \emph{non-degenerate} varieties
\cite[\S3]{inv.mot}. A smooth $G$-variety $X$ is non-degenerate
iff the action $G:T^{*}X$ is \emph{symplectically stable}, i.e.,
generic $G$-orbits in $\Im\Phi$ are closed in $\g^{*}$, see
Remark~\ref{symp.stab} below.

On the contrary, if $X$ is a generalized flag variety, then
$M=Q=G$ and $\Xx=X$. This is the opposite extremity in
Theorem~\ref{ref.loc.str}.

\begin{proof}[Proof of Theorem~\ref{ref.loc.str}]
In the notation of Theorem~\ref{gen.loc.str}, consider the
cotangent bundle over $\X$ and the conormal bundle $\U$ to the
foliation of $U$-orbits in~$\X$. We have
\begin{align*}
T^{*}\X\simeq\Ru{P}\times\Ru\p^{-}\times A\times\ab\times T^{*}C, \\
\U\simeq\Ru{P}\times\{0\}\times A\times\ab\times T^{*}C.
\end{align*}

\begin{claim}\label{Phi(U)}
$\Phi(\U)\subseteq\ab+\Ru\q$.
\end{claim}

Indeed, as $\Phi$ is equivariant and $\ab+\Ru\q$ is $P$-stable,
it suffices to consider~$\Phi(\U|_C)$. Take any $\alpha\in\U_x$,
$x\in C$, and let $\lambda$ be its projection to~$\ab$. We have
\begin{equation*}
(\Phi(\alpha),\xi)=\langle\alpha,\xi{x}\rangle=
\begin{cases}
      0,       & \xi\in\p_0; \\
(\lambda,\xi), & \xi\in\ab.
\end{cases}
\end{equation*}
Hence $\Phi(\alpha)\in\lambda+\Ru\p$. But $\alpha$ is fixed
by~$L_0$ whence
\begin{equation}\label{Phi(alpha)}
\Phi(\alpha)\in\lambda+\Ru\p^{L_0}\subseteq\lambda+\Ru\q,
\end{equation}
because $\Ru\p^{L_0}\cap\m=\Ru\p^L=0$.

\begin{claim}\label{sliding}
$\m_0{x}\subseteq\Ru\p{x},\ \forall x\in\X$.
\end{claim}

Otherwise there exists $\alpha\in\U_x$ that does not
vanish on~$\m_0{x}$. But then $\Phi(\alpha)\not\perp\m_0$, a
contradiction with Claim~\ref{Phi(U)}.

It follows from Claim~\ref{sliding} that
$\overline{\Ru{P}x}=Q_0x$, $\forall x\in\X$. Recall that the
orbits $\Ru{P}x=P_0x=Ux$ of points $x\in\X$ are parametrized
by~$\Z$. Hence the open set $\Xx=Q\X$ carries a foliation of
$Q_0$-orbits parametrized by~$\Z$ and $\U$ extends to the conormal
bundle of this foliation (denoted by the same letter). The
quotient space $\Xx/Q_0$ is $Q$-isomorphic to $\Z\simeq A\times C$
and $\U$ is the pullback of $T^{*}\Z\simeq A\times\ab\times
T^{*}C$. (Here $Q$ acts on $\Z$ through $A\simeq Q/Q_0$.)

We can lift $\Xx$ into $\U$ as the pullback of a $Q$-invariant
section $A\times\{\lambda\}\times C$ of $T^{*}\Z$
($\lambda\in\ab$). Formula~(\ref{Phi(alpha)}) yields a
commutative diagram of $Q$-equivariant maps
\begin{equation*}
\begin{CD}
T^{*}X @>{\Phi}>>      \g      \\
 @AAA                \bigcup   \\
  \Xx  @>{\phi}>> \lambda+\Ru\q.
\end{CD}
\end{equation*}

Now suppose $\lambda\in\ab^{\pr}$; then
$[\q,\lambda]=[\Ru\q,\lambda]=\Ru\q$. As
$Q\lambda=\Ru{Q}\lambda$ is closed in~$\g$
\cite[Satz~III.1.1-4]{IT}, we have $Q\lambda=\lambda+\Ru\q\simeq
Q/M$. It follows that $\Xx\simeq Q\itimes{M}\Zz$ is a
homogeneous fibering over $Q\lambda$ with fibre map $\phi$ and
fibre $\Zz=\phi^{-1}(\lambda)$.

\begin{claim}\label{flags}
$Mx\simeq M/M\cap P_0^{-},\ \forall x\in\Zz$.
\end{claim}

Indeed, the action $M\cap\Ru{P}:\Zz\cap\X$ is free and
$\overline{(M\cap\Ru{P})x}={Q_0x\cap\Zz}=M_0x$, $\forall
x\in\Zz\cap\X$. As $M_0\cap P$ normalizes $M\cap\Ru{P}$, we
obtain $(M\cap\Ru{P})x=(M_0\cap P)x$ and, without loss of
generality, $(M_0\cap P)_x=L_0$. Hence $(M_0)_x=M_0\cap P^{-}$,
the unique subgroup of $M_0$ containing $L_0$ and transversal to
$M\cap\Ru{P}$. Since $A\simeq M/M_0$ acts on
$\Zz/M_0\simeq\Xx/Q_0$ freely, we have $M_x=(M_0)_x=M\cap
P_0^{-}$.

Finally, replacing $C$ by the set of points in $\Zz$ with
stabilizer $M\cap P_0^{-}$ (which is a cross-section for
$Q$-orbits in~$\Xx$ by Claim~\ref{flags}) yields $\Zz\simeq
{(M/M\cap P_0^{-})\times C}$.
\end{proof}

\section{Geometry of cotangent bundle}\label{cotan}

Making use of the local structure theorem, we investigate here
the equivariant geometry of the cotangent bundle $T^{*}X$ of a
smooth $G$-variety~$X$.

Consider the conormal bundle to the foliation of generic
$\Ru{Q}$-orbits:
\begin{equation*}
\N=\{\alpha\in T^{*}_xX\mid x\in\Xx,\
\langle\alpha,\Ru\q{x}\rangle=0\}
\end{equation*}
By (\ref{X_1}) and (\ref{T*(G/H)}) we have
\begin{equation*}
\N\simeq\Ru{Q}\times T^{*}\Zz\simeq\Ru{Q}\times M\itimes{M\cap
P_0^{-}}(\ab+\m\cap\Ru\p^{-})\times T^{*}C.
\end{equation*}

\begin{lemma}\label{Phi(N)}
$\overline{\Phi(\N)}=\ab+\M+\Ru\q$, where
$\M=M(\m\cap\Ru\p^{-})=M(\m\cap\Ru\p)$ is the closure of a
Richardson nilpotent orbit in~$\m$.
\end{lemma}

\begin{proof}
Take any $\alpha\in\N_x$, $x\in C$. Let $\mu$ be the projection
of $\alpha$ to $N^{*}_x(\Ru{Q}C)\simeq
T^{*}_xMx\simeq\ab+\m\cap\Ru\p^{-}$, with the Jordan
decomposition $\mu=\sms\mu+\nil\mu$ ($\sms\mu\in\ab$,
$\nil\mu\in\m\cap\Ru\p^{-}$). We have
\begin{equation*}
(\Phi(\alpha),\xi)=\langle\alpha,\xi{x}\rangle=
\begin{cases}
    0,     & \xi\in\Ru\q; \\
(\mu,\xi), & \xi\in\m.
\end{cases}
\end{equation*}
Hence $\Phi(\alpha)\in\mu+\Ru\q$. It follows that
\begin{equation*}
\Phi(\N)=Q\Phi(\N|_C)\subseteq
Q(\ab+{\m\cap\Ru\p^{-}+\Ru\q})=\ab+\M+\Ru\q.
\end{equation*}

On the other hand, if $\sms\mu\in\ab^{\pr}$, then
$Q\mu=\sms\mu+M\nil\mu+\Ru\q$. Indeed, $Z(\mu)\subseteq
Z(\sms\mu)=M$ whence $[\Ru\q,\mu]=\Ru\q$ and
$\Ru{Q}\mu=\mu+\Ru\q$ by \cite[Satz~III.1.1-4]{IT}. Thus
$\Phi(\N)\supseteq\ab^{\pr}+\M+\Ru\q$ whence the assertion
on~$\overline{\Phi(\N)}$.

The equality $M(\m\cap\Ru\p^{-})=M(\m\cap\Ru\p)$ stems from a
well-known property of induced nilpotent orbits
\cite[7.1.3]{Ad}.
\end{proof}

\begin{corollary}\label{GN}
$\overline{G\N}=T^{*}X$.
\end{corollary}

\begin{proof}
Since $\N$ is $Q$-stable,
$\overline{G\N}=\overline{\Ru{Q^{-}}\N}$. Observe that
$\codim\N=\dim\Ru{Q}$. Hence it suffices to prove that
$\Ru{Q^{-}}\alpha\simeq\Ru{Q^{-}}$ is transversal to $\N$ for
some $\alpha\in\N$.

Choose $\alpha\in\N$ such that
$\mu=\Phi(\alpha)\in\ab^{\pr}+\M$. Then
$[\Ru\q^{-},\mu]=\Ru\q^{-}$ is transversal to~$\Phi(\N)$. The
assertion follows.
\end{proof}

\begin{corollary}[{\cite[5.4]{W_X}}]\label{Im(Phi)}
$\overline{\Im\Phi}=G(\ab+\M+\Ru\q)=\overline{G(\ab+\M)}
=G(\ab+\Ru\p)=G(\ab+\Ru\p^{-})$.
\end{corollary}

\begin{proof}
It remains only to note that $G(\ab+\M+\Ru\q)$ is closed,
because $\ab+\M+\Ru\q$ is stable under a parabolic $Q\subseteq
G$; similarly for $G(\ab+\Ru\p^{\pm})$.
\end{proof}

\begin{remark}\label{symp.stab}
In particular, generic $G$-orbits in $\Im\Phi$ are represented
by $\mu\in\ab+\M$. The orbit $G\mu$ is closed in $\g$ iff
$\mu=\sms\mu\in\ab$. Thus the action $G:T^{*}X$ is
symplectically stable iff $\M=0$ iff $M=L$.
\end{remark}

\begin{corollary}[{\cite[\S8]{W_X}}]\label{sgp}
The stabilizers of points in general position for the action
$G:T^{*}X$ are conjugate to the stabilizer of a point from the
open orbit for the action $M\cap P_0^{-}:\m\cap\Ru\p^{-}$.
\end{corollary}

\begin{proof}
Take $\alpha\in T^{*}X$ in general position. Without loss of
generality we may assume $\alpha\in\N$,
$\Phi(\alpha)=\mu\in\ab^{\pr}+\M$. Then $G_{\alpha}\subseteq
G_{\mu}\subseteq G_{\sms\mu}=M\subseteq Q$. However the
stabilizer in general position for $Q:\N$ coincides with that
for $M\cap P_0^{-}:\m\cap\Ru\p^{-}$.
\end{proof}

The image of the moment map contains the \emph{principal open
stratum} $(\Im\Phi)^{\pr}=G({\ab^{\pr}+\M})$. Let $T^{\pr}X$
denote its preimage in~$T^{*}X$, called the \emph{principal
stratum} of the cotangent bundle.

\begin{proposition}\label{T*X}
The principal stratum of $T^{*}X$ has the structure
\begin{equation}\label{T(pr)X}
T^{\pr}X\simeq G\itimes{N_X}\Sigma,
\end{equation}
where $\Sigma$ is the unique component of
$\Phi^{-1}({\ab^{\pr}+\M})$ intersecting $\N$
and $N_X$ is the stabilizer of $\Sigma$ in~$N(\ab)$.
\end{proposition}

\begin{proof}
It is easy to deduce from Definition~\ref{pr} that
\begin{equation*}
(\Im\Phi)^{\pr}\simeq G\itimes{N(\ab)}(\ab^{\pr}+N(\ab)\M)
\end{equation*}
and $N(\ab)\M$ is a union of Richardson orbit closures in $\m$
permuted transitively by the Weyl group $W(\ab)=N(\ab)/M$
of~$\ab$. Hence
\begin{equation*}
T^{\pr}X\simeq G\itimes{N(\ab)}\widetilde\Sigma,
\end{equation*}
where $\widetilde\Sigma=\Phi^{-1}(\ab^{\pr}+N(\ab)\M)$.
The fibre $\widetilde\Sigma$ is smooth and its components are
permuted transitively by~$W(\ab)$, because $T^{\pr}X$ is
irreducible. Thus (\ref{T(pr)X}) holds for any component
$\Sigma$ of~$\widetilde\Sigma$.

By Corollary~\ref{GN}, $\N$~intersects $T^{\pr}X$, and
Lemma~\ref{Phi(N)} implies that $\N^{\pr}=\N\cap
T^{\pr}X\subseteq\Ru{Q}\widetilde\Sigma\simeq
\Ru{Q}\times\widetilde\Sigma$. Since $\N$ is irreducible, it
intersects a unique component $\Sigma$ of~$\widetilde\Sigma$,
and $\Phi(\Sigma)=\ab^{\pr}+\M$ by Lemma~\ref{Phi(N)} again.
\end{proof}

\begin{remark}
The variety $\Sigma$ is a \emph{cross-section} of $T^{*}X$ in
the terminology of Guillemin--Sternberg \cite{symp} and Knop
\cite[5.4]{W(ham)}.
\end{remark}

\begin{remark}
The subgroup $W_X=N_X/M\subseteq W(\ab)$, i.e., the stabilizer
of $\Sigma$ in~$W(\ab)$, is nothing else but the \emph{little
Weyl group} of $X$ defined in~\cite{W_X}. To see this, consider
a morphism $\Psi:T^{\pr}X\to G\ab^{\pr}$ taking $\alpha$ to the
semisimple part of~$\Phi(\alpha)$. It splits as
\begin{equation}\label{Stein}
\begin{CD}
T^{\pr}X \simeq \itimes{N_X}\Sigma @>>> G\itimes{N_X}\ab^{\pr}
@>>> G\itimes{N(\ab)}\ab^{\pr} \simeq G\ab^{\pr}.
\end{CD}
\end{equation}
The right arrow is a finite morphism, and the left one has
irreducible generic fibres.

More precisely, restrict $\Psi$ to an open subset $\Sigma\cap
T^{*}\Xx=\Sigma\cap\N$ of~$\Sigma$. Since
$\N^{\pr}\simeq\Ru{Q}\times(\Sigma\cap\N)\simeq\Ru{Q}\times
T^{\pr}\Zz$, we have $M$-isomorphisms
\begin{equation}\label{Sigma*N}
\Sigma\cap\N\simeq T^{\pr}\Zz\simeq T^{*}(M_0/M_0\cap P^{-})
\times A\times\ab^{\pr}\times T^{*}C
\end{equation}
and $\Psi$ is just the projection to~$\ab^{\pr}$. So the fibres
are isomorphic to $A\times T^{*}(M_0C)$.

On the other hand, $\ab$~embeds into $\N$ as the conormal space
$N^{*}_x(Q_0C)$, $\forall x\in C$. Factoring (\ref{Stein}) by
$G$ we obtain a commutative diagram
\begin{equation*}
\begin{CD}
T^{\pr}X  @>>> G\itimes{N_X}\ab^{\pr} @>>> G\itimes{N(\ab)}\ab^{\pr} \\
  @AAA                  @VVV                         @VVV            \\
\ab^{\pr} @>>>     \ab^{\pr}/W_X      @>>>     \ab^{\pr}/W(\ab),
\end{CD}
\end{equation*}
so that the composite morphism $T^{\pr}X\to\ab^{\pr}/W_X$ has
irreducible generic fibres. Thus $W_X$ satisfies the definition
of \cite[\S6]{W_X}.
\end{remark}

\section{Horospheres and symplectic covering}\label{horosph}

\subsection{}\label{deg.hor}

Now we define a family of (possibly non-generic) horospheres
which is used below to polarize $T^{*}X$ and generalize
Theorem~\ref{cover(qaff)}.

Consider the set of $\Ru{Q}$-orbits in $\Xx$ (parametrized
by~$\Zz$) and the set $Y$ of their $G$-translates. One may say
that $Y$ is the set of generic orbits of unipotent subgroups
conjugate to~$\Ru{Q}$. On the other hand, $Y$~consists of
horospheres: every $\Ru{Q}$-orbit in $\Xx$ is an $M$-translate
of~$\Ru{Q}x$, $x\in\Zz$, such that $M_x\supset M\cap U$ whence
$\Ru{Q}x=Ux$.

\begin{proposition}\label{var(hor)}
The set $Y$ carries a structure of an algebraic $G$-variety such
that $Y\simeq G\itimes{Q}\Zz\simeq G/\Ru{Q}(M\cap P_0^{-})\times
C$ and $\dim Y=\dim X$. (Here $Q$ acts on $\Zz$ through
$M=Q/\Ru{Q}$.)
\end{proposition}

\begin{proof}
There is a natural map $G\itimes{Q}\Zz\to Y$, $g*x\mapsto
g\Ru{Q}x$. It suffices to verify that it is a bijection. The
assertion on dimensions is obvious since $\dim G\itimes{Q}\Zz=
\dim G/Q+\dim\Zz=\dim\Ru{Q}+\dim\Zz=\dim X$. The bijectivity
stems from the following
\let\qedsymbol\relax\end{proof}

\begin{lemma}
$g\Ru{Q}x=\Ru{Q}x'$ ($x,x'\in\Zz$) iff $g\in Q$.
\end{lemma}

\begin{proof}[Proof of the lemma]
As $\Zz=MC$ and $M$ normalizes~$\Ru{Q}$, we may assume without
loss of generality that $x,x'\in C$.

First suppose that $x=x'$. Let $S$ be the stabilizer of
$\Ru{Q}x$ in~$G$. Then $S=\Ru{Q}\cdot S_x\supseteq\Ru{Q}(M\cap
P_0^{-})$ contains a maximal unipotent subgroup of~$G$. The
structure of such groups is well known \cite[\S2]{W_X}: we have
$S=\Ru{S}\leftthreetimes\widetilde{S}$, where $\Ru{S}$ is the
unipotent radical of a parabolic subgroup $N(S)$ and
$\widetilde{S}$ is intermediate between the Levi subgroup of
$N(S)$ and its semisimple part. Furthermore,
$\Ru{S}\subseteq\Ru{Q}(M\cap\Ru{P}^{-})$ and
$\widetilde{S}\supseteq L_0$.

We have $\widetilde{S}=\widetilde{\Ru{Q}}\cdot\widetilde{S_x}$,
where $\widetilde{K}$ denotes the projection of $K\subseteq S$
to~$\widetilde{S}$. Hence $\widetilde{S}/\widetilde{S_x}$ is an
affine $\widetilde{S}$-variety with the transitive action of a
unipotent subgroup of~$\widetilde{S}$. This is possible only if
it is a point, i.e., $\widetilde{S}=\widetilde{S_x}$ whence
$S_x$ contains a conjugate of~$\widetilde{S}$. Therefore
$\widetilde{S}\subseteq S_y$ for some $y\in\Ru{Q}x$.

The subgroup $\widetilde{S}\cap M$ is reductive, contains~$L_0$,
and is contained in $S\cap M=M\cap P_0^{-}$. Hence
$\widetilde{S}\cap M=L_0$. If $\widetilde{S}\ne L_0$, then
$\widetilde{S}$ intersects~$\Ru{Q}$, a contradiction with
$\Ru{Q}\cap S_y=\emptyset$. It follows that $\widetilde{S}=L_0$,
$N(S)=\Ru{Q}(M\cap P^{-})$, and $S=\Ru{Q}(M\cap
P_0^{-})\subseteq Q$.

Finally, for arbitrary $x,x'\in C$ we have $g\in N(S)\subseteq
Q$.
\end{proof}

\subsection{}

Let $Y$ denote the variety of degenerate horospheres introduced
in~\ref{deg.hor}. Consider the polarized cotangent bundle with
respect to~$Y$:
\begin{equation*}
\widehat{T}^{*}X=\bigsqcup_{H\in Y}N^{*}H\simeq G\itimes{Q}\N.
\end{equation*}
(The latter isomorphism stems from Proposition~\ref{var(hor)}.)

\begin{remark}
The polarized cotangent bundle can be interpreted in a different
manner. Consider its principal open stratum
\begin{equation}\label{^T(pr)X}
\widehat{T}^{\pr}X=\hat{p}^{-1}(T^{\pr}X)\simeq
G\itimes{Q}\N^{\pr}\simeq G\itimes{M}(\Sigma\cap\N).
\end{equation}
By (\ref{T(pr)X}) the fibre product
\begin{equation}\label{fib.prod}
T^{\pr}X\ftimes{\ab^{\pr}/W_X}\ab^{\pr}\simeq G\itimes{M}\Sigma
\end{equation}
is birationally $G$-isomorphic to~$\widehat{T}^{*}X$. Thus we
generalize the definition of the polarized cotangent bundle in
\cite[\S3]{inv.mot}.
\end{remark}

Now we prove our second main result generalizing
Theorem~\ref{cover(qaff)}.

\begin{theorem}\label{cover}
There exists a $G$-equivariant symplectic rational Galois
covering $T^{*}Y\dasharrow T^{*}X$ with the Galois group~$W_X$.
\end{theorem}

\begin{proof}
We argue as in~\ref{hor}. It follows from (\ref{T(pr)X}) and
(\ref{^T(pr)X}) that $\hat{p}:\widehat{T}^{*}X\to T^{*}X$ is a
rational Galois covering with the Galois group~$W_X$. It remains
to prove that $\hat{q}:\widehat{T}^{*}X\to T^{*}Y$ is
birational.

The morphism $\hat{q}$ maps $N^{*}H$ to~$T^{*}_HY$, $\forall
H\in Y$. As $\hat{q}$ is equivariant, we may assume $H=\Ru{Q}x$,
$x\in C$. The action $\Ru{Q}:N^{*}H$ is free and any orbit
intersects $N^{*}_xH\simeq T^{*}_x\Zz$ in exactly one point,
i.e.,
\begin{equation*}
N^{*}H\simeq\Ru{Q}\times T^{*}_x\Zz\simeq
\Ru{Q}\times(\ab+\m\cap\Ru\p^{-})+T^{*}_xC.
\end{equation*}
On the other hand,
\begin{equation*}
T^{*}_HY\simeq (\g/(\Ru\q+\m\cap\p_0^{-}))^{*}\oplus T^{*}_xC
\simeq\ab+\m\cap\Ru\p^{-}+\Ru\q+T^{*}_xC
\end{equation*}
by Proposition~\ref{T*X}. It follows that the $\Ru{Q}$-action is
free on an open subset
\begin{equation*}
T^{\pr}_HY\simeq\ab^{\pr}+\m\cap\Ru\p^{-}+\Ru\q+T^{*}_xC
\end{equation*}
of~$T^{*}_HY$, and the orbits are just parallel planes with the
direction subspace~$\Ru\q$.

We already know that $\hat{p}$ is a rational covering, hence
$\hat{q}$ is a rational covering by Lemma~\ref{symp}. By
definition, $\hat{q}:N^{*}_xH\to N^{*}_HZ_x$ is a linear
isomorphism (where $Z_x\subset Y$ is the set of horospheres
containing~$x$). Hence $N^{*}_HZ_x$ intersects generic
$\Ru{Q}$-orbits in exactly one point, i.e., it is transversal to
$\Ru\q$ and projects onto
$\ab+\m\cap\Ru\p^{-}+T^{*}_xC\simeq T^{*}_x\Zz$ isomorphically.

More specifically, the composed map
\begin{equation*}
T^{*}_x\Zz\to N^{*}_xH\to N^{*}_HZ_x\to T^{*}_x\Zz
\end{equation*}
is identity. Indeed, the projection of $\alpha'\in N^{*}_xX$
or $\alpha''\in N^{*}_HZ_y$ to $T^{*}_x\Zz$ means the restriction
of $\alpha'$ to $T_x\Zz\subset T_xX$, or of $\alpha''$ to
$T_x\Zz\hookrightarrow T_HY$ (where $\Zz$ is regarded as a
subvariety of~$Y$), respectively. However, if we restrict the
construction of \ref{skew} to $\Zz\times\Zz\subset X\times Y$,
then the incidence variety $Z\subset X\times Y$ transforms to
$\diag\Zz$, so that $\alpha\in SN^{*}_{(x,x)}(\diag\Zz)$,
$\alpha'=\hat{p}(\alpha)$, $\alpha''=\hat{q}(\alpha)$ correspond
to one and the same covector in~$T^{*}_x\Zz$.

Finally, we conclude that
\begin{equation*}
\begin{CD}
N^{\pr}H\simeq\Ru{Q}\times(\ab^{\pr}+\m\cap\Ru\p^{-})+T^{*}_xC
@>{\hat{q}}>>T^{\pr}_HY
\end{CD}
\end{equation*}
is an isomorphism, which completes the proof.
\end{proof}

\section{Invariant collective motion}\label{inv.mot}

\subsection{}

The pullbacks of functions on $\g^{*}$ along $\Phi$ are called
\emph{collective functions} on~$T^{*}X$. Their skew gradients
generate the tangent spaces to $G$-orbits at every point. Hence
collective functions are in involution with $G$-invariant
functions on $T^{*}X$, i.e., they serve as simultaneous
integrals for all $G$-invariant Hamiltonian dynamical systems
on~$T^{*}X$.

Invariant collective functions, i.e., pullbacks of invariant
functions on~$\g^{*}$, have the property that their skew
gradients are both tangent and skew orthogonal to $G$-orbits,
and even generate the kernel of the symplectic form on
$\g{\alpha}$ for $\alpha\in T^{*}X$ in general position
\cite{symp}. Since invariant collective functions are in
involution, their skew gradients generate an Abelian flow of
$G$-equivariant symplectomorphisms of $T^{*}X$ preserving
$G$-orbits, which is called the \emph{invariant collective
motion}.

Restricted to any orbit $G\alpha\subset T^{*}X$, the invariant
collective motion gives rise to a connected Abelian subgroup of
$G$-automorphisms $A_{\alpha}\subseteq
N(G_{\alpha})/G_{\alpha}$. It is known \cite{symp} that
$A_{\alpha}\simeq\conn{(G_{\Phi(\alpha)}/G_{\alpha})}$ for
$\alpha$ in general position, cf.\ Remark~\ref{gen.flats} below.

A tempting problem is to integrate the invariant collective
motion, i.e., to find an algebraic group of $G$-equivariant
symplectomorphisms of $T^{*}X$ whose restriction to every orbit
$G\alpha$ coincides with~$A_{\alpha}$. However the problem has
generally no solution in this formulation, as the following
example shows.

\begin{example}
Let $X=G$ and $G$ act on $X$ by left translations. Then
$T^{*}X\simeq G\times\g^{*}$ is a trivial bundle with the
$G$-action by left translations of the first factor, and the
moment map $\Phi:G\times\g^{*}\to\g^{*}$ is just the coadjoint
action map.

Invariant collective functions on $T^{*}X$ are of the form
$F(g,\mu)=f(\mu)$, where $f$ is a $G$-invariant function
on~$\g^{*}$. The differential $dF$ vanishes along
$G\times\{\mu\}$ and equals $df$ on $\{e\}\times\g^{*}$. It
follows that the skew gradient of $F$ at $(e,\mu)$ is
$(d_{\mu}f,0)$. (Note that $d_{\mu}f\in\g^{**}\simeq\g$.)
If $\mu$ is a regular point (i.e., $\dim G\mu=\max$), then
the $d_{\mu}f$ span the annihilator $\g_{\mu}$
of~$\g\mu$, whence $\ab_{(g,\mu)}\simeq\g_{\mu}$ and
$A_{(g,\mu)}\simeq\conn{(G_{\mu})}$.

On the principal stratum $T^{\pr}X\simeq G\times\g^{\pr}$ (where
$\g^{\pr}$ is the set of all regular semisimple elements in
$\g\simeq\g^{*}$)  all groups $A_{(g,\mu)}$ are maximal tori
of~$G$. Hence the hypothetical automorphism group integrating
the invariant collective motion must be a torus. However for any
regular $\mu\notin\g^{\pr}$ the group $A_{(g,\mu)}$
contains unipotent elements. Thus one might hope to integrate
the invariant collective motion only on a proper open
subset~$T^{\pr}X$.

But even there it is not possible, because the family of tori
$A_{(g,\mu)}=G_{\mu}$ cannot be trivialized globally on
$G\times\g^{\pr}$. Indeed, the family of~$G_{\mu}$,
$\mu\in\g^{\pr}\simeq G\itimes{N(\ab)}\ab^{\pr}$, has the
non-trivial monodromy group $W=W(\ab)$, the Weyl group of~$\g$.
(Here $A$ is a maximal torus of~$G$.) Only if we unfold
$T^{\pr}X$ by taking a Galois covering of $\g^{\pr}$ by
$G\itimes{A}\ab^{\pr}\simeq G/A\times\ab^{\pr}$ with the Galois
group $W=N(\ab)/A$, then the family of $A_{(g,\mu)}$ lifts to a
trivial family of tori, so that we can integrate the invariant
collective motion on the covering space.
\end{example}

This example suggests the following reformulation of the
problem: to integrate the invariant collective motion on a
suitable {\'e}tale covering of an open subset in~$T^{*}X$. This
problem was solved by Knop \cite{inv.mot} for
non-degenerate~$X$. Here we consider the general case.

Recall from (\ref{T(pr)X}) and (\ref{^T(pr)X}) that the
polarization map $\hat{p}:\widehat{T}^{\pr}X\to T^{*}X$ is a
$G$-equivariant {\'e}tale Galois covering of an open subset
$G\N^{\pr}$ of~$T^{\pr}X$. The symplectic structure on $T^{*}X$
lifts to $\widehat{T}^{\pr}X$, so that $G:\widehat{T}^{\pr}X$ is
a Hamiltonian action with the moment map
$\widehat\Phi=\Phi\circ\hat{p}$. Hence the invariant collective
motion on $T^{*}X$ lifts to the invariant collective motion
on~$\widehat{T}^{\pr}X$.

We have a commutative diagram:
\begin{equation}\label{pol&cotan}
\begin{CD}
\llap{$\widehat{T}^{\pr}X\underset{\text{open}}\hookrightarrow\;$}
G\itimes{M}\Sigma @>\hat{p}>\text{{\'e}tale}> G\itimes{N_X}\Sigma
&\simeq& T^{\pr}X \\
       @VVV                                       @VV{\Phi}V @.\\
G\itimes{M}(\ab^{\pr}+\M) @>>\text{{\'e}tale}> G(\ab^{\pr}+\M)
&\;=\;& (\Im\Phi)^{\pr} \\
       @VVV                                          @VVV    @.\\
    \ab^{\pr}        @>>\text{{\'e}tale}>      \ab^{\pr}/W(\ab).
\end{CD}
\end{equation}
Note that the collective invariant functions on $T^{\pr}X$ are
exactly those pulled back from~$\ab^{\pr}/W(\ab)$. Consider the
composed map
\begin{equation*}
\Pi:\widehat{T}^{\pr}X \longrightarrow \ab^{\pr}, \quad
\Pi(g*\alpha)=\sms{\Phi(\alpha)}.
\end{equation*}

\begin{theorem}
There is a Hamiltonian $G$-equivariant $A$-action on
$\widehat{T}^{\pr}X$ with the moment map $\Pi$ which integrates
the invariant collective motion.
\end{theorem}

\begin{remark}
Instead of~$\widehat{T}^{\pr}X$ we may consider the fibre
product~(\ref{fib.prod}), which covers the whole $T^{\pr}X$ and
contains $\widehat{T}^{\pr}X$ as an open subset.
\end{remark}

\begin{proof}
The torus $Z(M)$ acts on $\Sigma$ with kernel
${Z(M)\cap L_0}$. Indeed,
\begin{equation*}
\Sigma\cap\N\simeq\N^{\pr}|_{\Zz}\simeq T^{\pr}(M/M\cap P_0^{-})
\times T^{*}C
\end{equation*}
by~(\ref{Sigma*N}), and the kernel of $Z(M):M/M\cap P_0^{-}$ is
exactly~$Z(M)\cap L_0$. This yields an action of
$A=Z(M)/(Z(M)\cap L_0)$ on $\Sigma$ which is free
on~$\Sigma\cap\N$. The $A$-action commutes with $M$ and
immediately extends to the whole $\widehat{T}^{\pr}X\simeq
G\itimes{M}(\Sigma\cap\N)$ by $G$-equivariance:
$a(g*\alpha)=g*z\alpha$, where $z\in Z(M)$ represents $a\in A$.

Now we prove that $\Pi$ is the moment map for this action.
Obviously $\Pi$ is $A$-invariant. For any $\xi\in\ab$ consider
two functions $\Pi^{*}\xi$ and $\widehat\Phi^{*}\xi$
on~$\widehat{T}^{\pr}X$. (We think of $\xi$ as of a linear
function on $\ab^{*}\simeq\ab$, or on $\g^{*}\simeq\g$,
respectively.)

\begin{Claim}
$d_{\alpha}\Pi^{*}\xi=d_{\alpha}\widehat\Phi^{*}\xi,\
\forall\xi\in\ab,\ \alpha\in\Sigma$
\end{Claim}

Indeed,
$T_{\alpha}\widehat{T}^{\pr}X=\g{\alpha}+T_{\alpha}\Sigma$. The
differentials coincide on $T_{\alpha}\Sigma$, because the
functions coincide on $\Sigma$:
\begin{equation*}
\Pi^{*}\xi(\alpha)=(\xi,\sms{\Phi(\alpha)})=
(\xi,\Phi(\alpha))=\widehat\Phi^{*}\xi(\alpha),
\qquad\forall\alpha\in\Sigma.
\end{equation*}
On the other hand, $d(\Pi^{*}\xi)$~vanishes on $\g{\alpha}$ for
$\Pi^{*}\xi$ is $G$-invariant, and
\begin{equation*}
d\widehat\Phi^{*}\xi(\g\alpha)=(\xi,d\widehat\Phi(\g\alpha))=
(\xi,[\g,\Phi(\alpha)])=
(\xi,[\m_0,\nil{\Phi(\alpha)}]+\Ru\q+\Ru\q^{-})=0.
\end{equation*}

It follows from the claim that the skew gradient of $\Pi^{*}\xi$
at any $\alpha\in\Sigma$ coincides with that of
$\widehat\Phi^{*}\xi$, i.e., with~$\xi\alpha$. By
$G$-invariance, we conclude that the skew gradient of
$\Pi^{*}\xi$ on $\widehat{T}^{\pr}X$ is the velocity field of
$\xi$ with respect to the above $A$-action. Thus the $A$-action
is symplectic and $\Pi$ is the moment map.

Finally, since the horizontal arrows in (\ref{pol&cotan}) are
{\'e}tale maps, the skew gradients of $\Pi^{*}\xi$ ($\xi\in\ab$)
span the same subspace in $T_{\alpha}\widehat{T}^{\pr}X$ as the
skew gradients of invariant collective functions. It follows
that the $A$-action integrates the invariant collective motion.
\end{proof}

Now we have three symplectic actions on~$\widehat{T}^{\pr}X$:
of~$G$, of~$A$ (integration of the invariant collective motion),
and of~$W_X$ (the Galois group). They patch together in the
following picture.

\begin{Corollary}
There is a Hamiltonian action $G\times(W_X\rightthreetimes
A):\widehat{T}^{\pr}X$ with the moment map
$(\widehat\Phi,\Pi):\widehat{T}^{\pr}X\to\g^{*}\oplus\ab^{*}$.
\end{Corollary}

\begin{proof}
The group $N_X$ acts on $\Sigma$ and on $Z(M)$ by conjugation.
Hence $N_X$ preserves the kernel $Z(M)\cap L_0$ of
$Z(M):\Sigma$, i.e., $W_X$~acts on $A=Z(M)/(Z(M)\cap L_0)$. The
actions of $W_X$ and $A$ on $\widehat{T}^{\pr}X$ patch together
into the $(W_X\rightthreetimes A)$-action, as the following
calculation shows: let $n\in N_X$, $z\in Z(M)$ represent $w\in
W_X$, $a\in A$, respectively; then
\begin{multline*}
w\cdot a\cdot(g*\alpha)=w(g*z\alpha)=gn^{-1}*nz\alpha=\\
gn^{-1}*(nzn^{-1})n\alpha=(wa)\cdot(gn^{-1}*n\alpha)
=(wa)\cdot w(g*\alpha).
\end{multline*}
It remains to note that $\Pi$ is $W_X$-equivariant.
\end{proof}

\subsection{}

The orbits of the invariant collective motion in $T^{\pr}X$
(or~$\widehat{T}^{\pr}X$) can be defined intrinsically as
$Z(M(\alpha))\cdot\alpha$, where $M(\alpha)=Z(\sms\mu)$ for
${\mu=\Phi(\alpha)}$ (or $\mu=\widehat\Phi(\alpha)$). Note that
the $A$-action on $\widehat{T}^{\pr}X$ is free, hence
$Z(M(\alpha))\,\alpha\simeq A$,
$\forall\alpha\in\widehat{T}^{\pr}X$.

The projections of these orbits to $X$ are called \emph{flats}.
Namely, a flat through $x\in X$ is $F_{\alpha}=Z(M(\alpha))\,x$,
where $\alpha$ is any (polarized) covector over $x$ with
$\Phi(\alpha)\in(\Im\Phi)^{\pr}$
(or $\widehat\Phi(\alpha)\in(\Im\Phi)^{\pr}$, respectively).
Clearly, $F_{\alpha}$ does not depend on the polarization
of~$\alpha$.

It is easy to see that $F_{\alpha}\simeq A$,
$\forall\alpha\in\widehat{T}^{\pr}X$. Indeed, without loss of
generality we may assume $\alpha\in\Sigma\cap\N$, $x\in\Xx$.
Then $M(\alpha)=M$, and $Z(M)_x\subseteq L_0$ by
Theorem~\ref{ref.loc.str}. Hence $Z(M)_x=Z(M)_{\alpha}=Z(M)\cap
L_0$, which implies the claim.

\begin{remark}\label{gen.flats}
If $\alpha\in\widehat{T}^{\pr}X$ is in general position, namely,
the $G$-orbit of $\mu=\widehat\Phi(\alpha)$ has maximal
dimension in~$\Im\Phi$, then the orbit of the invariant
collective motion coincides with $\conn{(G_{\mu})}\alpha$ and
$F_{\alpha}=\conn{(G_{\mu})}x$, cf.~\cite{symp},
\cite[\S4]{inv.mot}. Indeed, assuming $\alpha\in\Sigma\cap\N$,
$\mu\in\ab^{\pr}+\m\cap\Ru\p^{-}$, we have $G_{\mu}=M_{\nil\mu}$
and $M\nil\mu$ is open in~$\M$. Then
$\conn{(M_{\nil\mu})}\subseteq M\cap P^{-}$ and $(M\cap
P^{-})_{\nil\mu}=Z(M)\cdot(M\cap P_0^{-})_{\nil\mu}=Z(M)\cdot
G_{\alpha}$.
\end{remark}

Now let $X\hookrightarrow\overline{X}$ be a $G$-equivariant open
embedding. The closures of flats in $\overline{X}$ are (possibly
non-normal) toric varieties. Rigidity of tori implies that the
closures of generic flats are isomorphic. It is an interesting
problem to describe the closure of a generic flat. This has
important applications in the equivariant embedding theory,
see~\cite{inv.mot}. The problem was solved by Knop for
non-degenerate $X$~\cite{inv.mot}, and it would be desirable to
extend his solution to arbitrary~$X$.

\end{document}